\documentclass[onefignum,onetabnum]{siamart171218}

\usepackage{lipsum}
\usepackage{amsfonts}
\usepackage{graphicx}
\usepackage{epstopdf}
\usepackage{algorithmic}
\usepackage{subcaption}
\ifpdf
  \DeclareGraphicsExtensions{.eps,.pdf,.png,.jpg}
\else
  \DeclareGraphicsExtensions{.eps}
\fi

\newsiamremark{remark}{Remark}
\newsiamremark{hypothesis}{Hypothesis}
\crefname{hypothesis}{Hypothesis}{Hypotheses}
\newsiamthm{claim}{Claim}

\headers{Efficient Explicit Taylor ODE Integrators with Symbolic-Numeric Computing}{S. Tan, O. Smith and C. Rackauckas}

\title{Efficient Explicit Taylor ODE Integrators with Symbolic-Numeric Computing}

\author{
  Songchen Tan\thanks{Department of Mathematics and Center for Computational Science and Engineering, Massachusetts Institute of Technology, Cambridge, MA 
  (\email{songchen@mit.edu}).}
  \and
  Oscar Smith\thanks{JuliaHub, Cambridge, MA (\email{oscar.smith@juliahub.com}).}
  \and
  Christopher Rackauckas\thanks{Department of Mathematics, Massachusetts Institute of Technology, Cambridge MA 
  (\email{crackauc@mit.edu}).}
}

\usepackage{amsopn}

\makeatletter
\newcommand*{\addFileDependency}[1]{% argument=file name and extension
  \typeout{(#1)}% latexmk will find this if $recorder=0 (however, in that case, it will ignore #1 if it is a .aux or .pdf file etc and it exists! if it doesn't exist, it will appear in the list of dependents regardless)
  \@addtofilelist{#1}% if you want it to appear in \listfiles, not really necessary and latexmk doesn't use this
  \IfFileExists{#1}{}{\typeout{No file #1.}}% latexmk will find this message if #1 doesn't exist (yet)
}
\makeatother

\ifpdf
\hypersetup{
  pdftitle={An Example Article},
  pdfauthor={D. Doe, P. T. Frank, and J. E. Smith}
}
\fi

\begin{document}

\maketitle

% REQUIRED
\begin{abstract}
  Taylor series methods show a newfound promise for the solution of non-stiff ordinary differential equations (ODEs) given the rise of new compiler-enhanced techniques for calculating high order derivatives. In this paper we detail a new Julia-based implementation that has two important techniques: (1) a general purpose higher-order automatic differentiation engine for derivative evaluation with low overhead; (2) a combined symbolic-numeric approach to generate code for recursively computing the Taylor polynomial of the ODE solution. We demonstrate that the resulting software's compiler-based tooling is transparent to the user, requiring no changes from interfaces required to use standard explicit Runge-Kutta methods, while achieving better run time performance. In addition, we also developed a comprehensive adaptive time and order algorithm that uses different step size and polynomial degree across the integration period, which makes this implementation more efficient and versatile in a broad range of dynamics. We show that for codes compatible with compiler transformations, these integrators are more efficient and robust than the traditionally used explicit Runge-Kutta methods.
\end{abstract}

% REQUIRED

% % REQUIRED
% \begin{AMS}
%   68Q25, 68R10, 68U05
% \end{AMS}

\section{Introduction}

The numerical integration of ordinary differential equations (ODEs) remains a cornerstone of computational science and engineering, with applications ranging from celestial mechanics and astrodynamics~\cite{battin_introduction_1999}, fluid dynamics~\cite{ferziger_computational_2002}, and control systems~\cite{astrom_feedback_2010} to modern machine learning models such as neural ODEs~\cite{chen_neural_2018}. 
Among diverse classes of integrators for non-stiff problems, methods based on Taylor polynomials (also known as Taylor methods)~\cite{corliss_solving_1982} offer an intriguing alternative to the more widely employed Runge–Kutta schemes~\cite{noauthor_solving_1993,shampine_solving_2002,shampine_matlab_1997,vladimirescu_spice_1994}. 
In an explicit Taylor method, one computes a Taylor polynomial of the solution about the current time and then uses that polynomial to advance the solution for a step size. Formally, for an initial value problem 
$$
u'=f(u,t);\quad u(t_0)=u_0
$$
with $f:(\mathbb R^d\times\mathbb R)\to\mathbb R^d$ and $u_0\in\mathbb R^d$, assuming the discretized times are $t_n$ for $n\in\mathbb N$, corresponding solution values are $u_n=u(t_n)$, derivatives are $u_n^{(k)}=(\mathrm d^ku/\mathrm dt^k)(t_n)$. At each step, one computes a polynomial of degree $p$ and step size $h$:
$$
u_{n+1}=u_n+u^{(1)}_nh+\frac12u^{(2)}_nh^2+\cdots+\frac{1}{p!}u^{(p)}_nh^p
$$
For the rest of this paper we will omit the word ``explicit'' and refer to this as the ``Taylor method''.

The appeal of Taylor methods lies chiefly in two features. Firstly, for sufficiently smooth functions $f$ one can systematically increase the convergence order of the method by computing more terms in the polynomial (via automatic differentiation). For problems requiring very high precision (e.g., long‐term orbit propagation, sensitivity analysis), high‐order Taylor methods can outperform standard methods~\cite{jorba_software_2005,baeza_stability_2018}. Secondly, the polynomial interpolation of solution $u$ comes for free in the Taylor approach, therefore one can readily compute $u(t)$ for arbitrary $t\in(t_n,t_{n+1})$ and construct a dense output~\cite{abad_algorithm_2012}.

Despite these advantages, Taylor methods have been less popular because of the complexity of evaluating higher-order derivatives, while the most widely used explicit Runge-Kutta methods only require that the user supplies the function $f$ to the solver. For evaluating higher-order derivatives, because numerical differentiation is too slow and inaccurate, Taylor methods need some way to compute higher order derivatives exactly, which either requires user intervention to derive such functions or some automatic differentiation via algorithmic, symbolic, or compiler means. The chief challenge is to construct an automatic differentiation engine that can compute arbitrary high order derivatives of a broad range of functions with small overhead~\cite{chang_atomft_1994}. In addition, computing the Taylor polynomial of an ODE solution is a recursive process in nature~\cite{griewank_evaluating_2008}, which must be combined in a clever way with the aforementioned higher-order automatic differentiation engine. If, instead, the automatic differentiation engine is regarded as a separate ``black box'' from the ODE integrator, there will be a lot of redundancies. To see this, one can imagine the procedure for $p=3$:
\begin{enumerate}
    \item Obtain $u_n^{(1)}$ using the relation $u_n^{(1)}=f(u(t),t)|_{t=t_n}$;
    \item Obtain $u_n^{(2)}$ using the relation $u_n^{(2)}=(\mathrm d/\mathrm dt)(f(u(t),t))|_{t=t_n}$, which needs $u_n^{(1)}$ to be already available;
    \item Obtain $u_n^{(3)}$ using the relation $u_n^{(3)}=(\mathrm d^2/\mathrm dt^2)(f(u(t),t))|_{t=t_n}$, which needs $u_n^{(1)},u_n^{(2)}$ to be already available;
\end{enumerate}
Note that here $\mathrm d/\mathrm dt$ means the total derivative
$$
\frac{\mathrm d}{\mathrm dt}f(u(t),t)=\frac{\partial f}{\partial u}\frac{\mathrm du}{\mathrm dt}+\frac{\partial f}{\partial t}
$$
Therefore, one can only get one additional term with each call to the higher-order AD engine. Assuming that doing order-$k$ AD on $f$ needs $O(k^2)$ times the operations in $f$ (which is a reasonable assumption for most implementation of Taylor-mode AD), the total complexity should be $O(1^2+2^2+\cdots+p^2)=O(p^3/3)$ which is prohibitively expensive~\cite{griewank_evaluating_2008}.

To overcome this problem, previous implementations such as ATOMFT~\cite{berryman_atomft_1988,chang_atomft_1994}, \verb|taylor| by Jorba \emph{et al}~\cite{jorba_software_2005} and TIDES~\cite{abad_algorithm_2012} must develop their special-purpose automatic differentiation engine specifically for the application of solving ODEs, which can efficiently handle the recursive computation of Taylor polynomial. These approaches are asymptotically efficient since it is proven ~\cite{griewank_evaluating_2008} that the complexity could be reduced to $O(p^2)$. However, using languages like C and FORTRAN, they often need a separate build step to transform user-defined functions to the actual code computing Taylor polynomials of ODE solutions, either via a hand-written parser~\cite{chang_atomft_1994} or another computer algebra system like Mathematica~\cite{abad_algorithm_2012}. The build step imposes significant limitations on how the user could write their ODE function (e.g. they can only use a small set of elementary functions), and makes these implementations hard to incorporate in modern computing environments where scientific and engineering users often call the methods from a dynamical language interface. Within the Julia community, TaylorIntegration.jl~\cite{perez_perezhztaylorintegrationjl_2025} uses a general purpose higher-order automatic differentiation engine~\cite{benet_taylorseriesjl_2019}, but relies on \emph{ad hoc} macros to analyze user-defined functions and partially remove the redundancy in computing Taylor polynomials. This implementation reduces the barrier to use Taylor methods, but still requires users to write ODE functions in a pretty restrictive way such that it could be parsed by the specific macro. We will compare this approach with our approach in more detail in the appendix.

In addition to these transformation-based methods, Griewank \emph{et al}~\cite{griewank_evaluating_2008} also introduced a coefficient doubling method that could also reduce the complexity to $O(p^2)$, but at the complexity of concretely storing some Jacobian matrices $A_i\equiv(\partial u_i/\partial u_0)$. In the context of solving large-scale ODE systems, this is not appealing since the space complexity on dimension $d$ could be $O(d^2)$ instead of $O(d)$. The JAX implementation~\cite{bettencourt_taylor-mode_2022} of this algorithm tried to avoid storing Jacobians by computing more Jacobian-vector products, but their performance result didn't show significant improvement compared to recursively computing coefficients.

Nevertheless, recent advances in higher-order automatic differentiation engines~\cite{tan_higher-order_2023} and symbolic-numeric computational frameworks~\cite{gowda_high-performance_2021} have given Taylor methods new possibilities. In this paper, we utilize these new techniques in developing a new implementation of the Taylor methods of arbitrary order. Our key contributions are:

\begin{itemize}
    \item Non-intrusively composed a general-purpose computer algebra system~\cite{gowda_high-performance_2021} with a general-purpose higher-order automatic differentiation engine~\cite{tan_taylordiffjl_2023} to establish an efficient, low-overhead framework for Taylor methods in Julia;
    \item Symbolically removed redundant computation and synthesized code for computing the Taylor polynomial of ODE solution. Compared with previous approaches, this provides a more seamless user experience since the code transformations happen at the user's first call to the solver using just-in-time compilation, and does not require the user to leave the usual Julia programming environment;
    \item Introduced and systematically evaluated an mechanism to use different step size and polynomial degree across different integration regions. In problems where the convergence radius of Taylor series of ODE functions might vary a lot, the adaptive polynomial degree solvers could outperform their fixed-degree counterparts.
\end{itemize}
The remainder of the paper is organized as follows. In section 2 we introduce the Taylor-mode AD and the symbolic-numeric approach to compute the Taylor polynomial of ODE solutions, as well as a mechanism to control the step size and polynomial degree across the integration time period. Section 3 shows numerical experiments with focus on comparing our integrator to state-of-the-art Runge-Kutta algorithms, and section 4 offers concluding remarks and directions for future research.

\section{Theory}\label{sec:theory}

\subsection{Taylor-mode automatic differentiation for higher-order directional derivatives}
Let $U$ be an open subset of $\mathbb R^{d_1}$, and let $f:U\to\mathbb R^{d_2}$ be a sufficiently smooth function. The derivative of $f$ at a point $x\in\mathbb R^{d_1}$ is a linear operator $Df(x):\mathbb R^{d_1}\to\mathbb R^{d_2}$ such that it maps a vector $v\in\mathbb R^n$ to $Df(x)[v]$, which is the directional derivative of $f$ at $x$ in the direction of $v$, also known as the Jacobian-vector product.

Similarly, for $p\in\mathbb Z^+>1$, the derivative of $p$-th order of $f$ at a point $x$ is a multilinear operator $D^pf(x):(\mathbb R^{d_1},\cdots,\mathbb R^{d_1})\to\mathbb R^{d_2}$ such that it maps a tuple of vectors $(v_1,\ldots,v_p)$ to $D^pf(x)[v_1,\cdots,v_p]$, the directional derivative of $f$ at $x$ in the direction of $(v_1,\ldots,v_p)$.

First-order forward-mode AD can be viewed as an algorithm to propagate directional derivative information through the composition of functions~\cite{revels_forward-mode_2016}. For example, let $f(x)=g(h(x))$ be a composite function for $\mathbb R^{d_1}\to\mathbb R^{d_2}$ and $x\in\mathbb R^{d_1}$ be a specific point in its domain. Suppose that we already have $h_0=h(x)$ representing the primal output of $h$, and $h_1=Dh(x)[v]$ representing the perturbation of $h$ along some direction $v\in\mathbb R^{d_1}$. Now for the composite function, we want to know what is the perturbation of $f$ along direction $v$; in other words, we want to know $f_1=Df(x)[v]$ in addition to $f_0=f(x)$. The answer to that is simply the chain rule:
\begin{itemize}
  \item $f_0 = g(h_0)$
  \item $f_1 = Dg(h_0)[h_1]$
\end{itemize}
Similarly, the Taylor-mode AD can be viewed as an algorithm to propagate higher-order directional derivative information through compositions of functions. Suppose that we already have a Taylor bundle for $h$ at $x$ along $v$:
$$
(h_0,h_1,\ldots,h_p) = (h(x),Dh(x)[v],D^2h(x)[v,v],\ldots,D^ph(x)[v,\ldots,v])
$$
and we want to know
$$
(f_0,f_1,\ldots,f_p) = (f(x),Df(x)[v],D^2f(x)[v,v],\ldots,D^pf(x)[v,\ldots,v])
$$
the answer to that is, in turn, the Fa\`a di Bruno's formula:
\begin{itemize}
  \item $f_0 = g(h_0)$
  \item $f_1 = Dg(h_0)[h_1]$
  \item $f_2 = D^2g(h_0)[h_1,h_1] + Dg(h_0)[h_2]$
  \item $\cdots$
\end{itemize}
In our practical implementation of Taylor-mode AD, a pushforward rule is defined for every ``simple'' function, so that derivatives of any complicated functions that are composed of these simple functions can be computed automatically, either via operator-overloading or source-code transformation techniques. Our previous work ~\cite{tan_higher-order_2023} has shown that for any functions composed of elementary functions and arbitrary control flow, the $p$-th order directional derivative is at most $O(p^2)$ times more expensive than computing the function itself, in contrast to naively nesting first-order forward-mode AD which could be $O(2^p)$ times expensive~\cite{bettencourt_taylor-mode_2022}. We provided a reference implementation of this method in the Julia language~\cite{bezanson_julia_2017} using its multiple-dispatch mechanism~\cite{bezanson_fast_2018}, TaylorDiff.jl~\cite{tan_taylordiffjl_2023}, featuring an automatic generation~\cite{gowda_high-performance_2021} of higher-order pushforward rules from first-order pushforward rules defined in ChainRules.jl~\cite{white_juliadiffchainrulesjl_2023}.

Below, we will demonstrate how to use Taylor-mode AD to implement ODE solvers.

\subsection{Computing Taylor polynomial of ODE solutions}

In the ODE $u'=f(u)$, the degree-$p$ Taylor polynomial of ODE solution is defined as
$$
T_n(t)=u_n+u_n^{(1)}(t-t_n)+\frac12u_n^{(2)}(t-t_n)^2+\cdots+\frac1{p!}u_n^{(p)}(t-t_n)^p+O((t-t_n)^{p+1})
$$
such that the integration step $u_{n+1}:=T_n(t_{n+1})$ has convergence order $p$. If we expand both $u(t)$ and $f(u(t))$ as Taylor series of $t$ near some $t_n$, then we get
$$
u_n^{(k+1)}=\frac{\mathrm d^{k+1}u}{\mathrm dt^{k+1}}(t_n)=\left[\frac{\mathrm d^k}{\mathrm dt^k}f(u(t),t)\right]\bigg |_{t=t_n},\quad\forall k\in\{0,\cdots,p-1\}
$$
Since $f$ is a function of $u$ and $t$, evaluating the total derivatives through the chain rules (i.e. Taylor-mode automatic differentiation) would require $u_n,u_n^1,\cdots,u_n^k$. Therefore, a naive approach is to do automatic differentiation on $f$ for multiple times, each time one order higher to obtain $(\mathrm d^k/\mathrm dt^k)f(u(t),t)|_{t=t_n}$ and then convert to $u_n^{(k+1)}$. However, this will result in a lot of redundant computation, because when $(\mathrm d^k/\mathrm dt^k)f(u(t),t)|_{t=t_n}$ is computed, $(\mathrm d^i/\mathrm dt^i)f(u(t),t)|_{t=t_n}$ for all $i<k$ are also computed again.

\begin{figure}[!ht]
    \centering
    \includegraphics[width=0.7\linewidth]{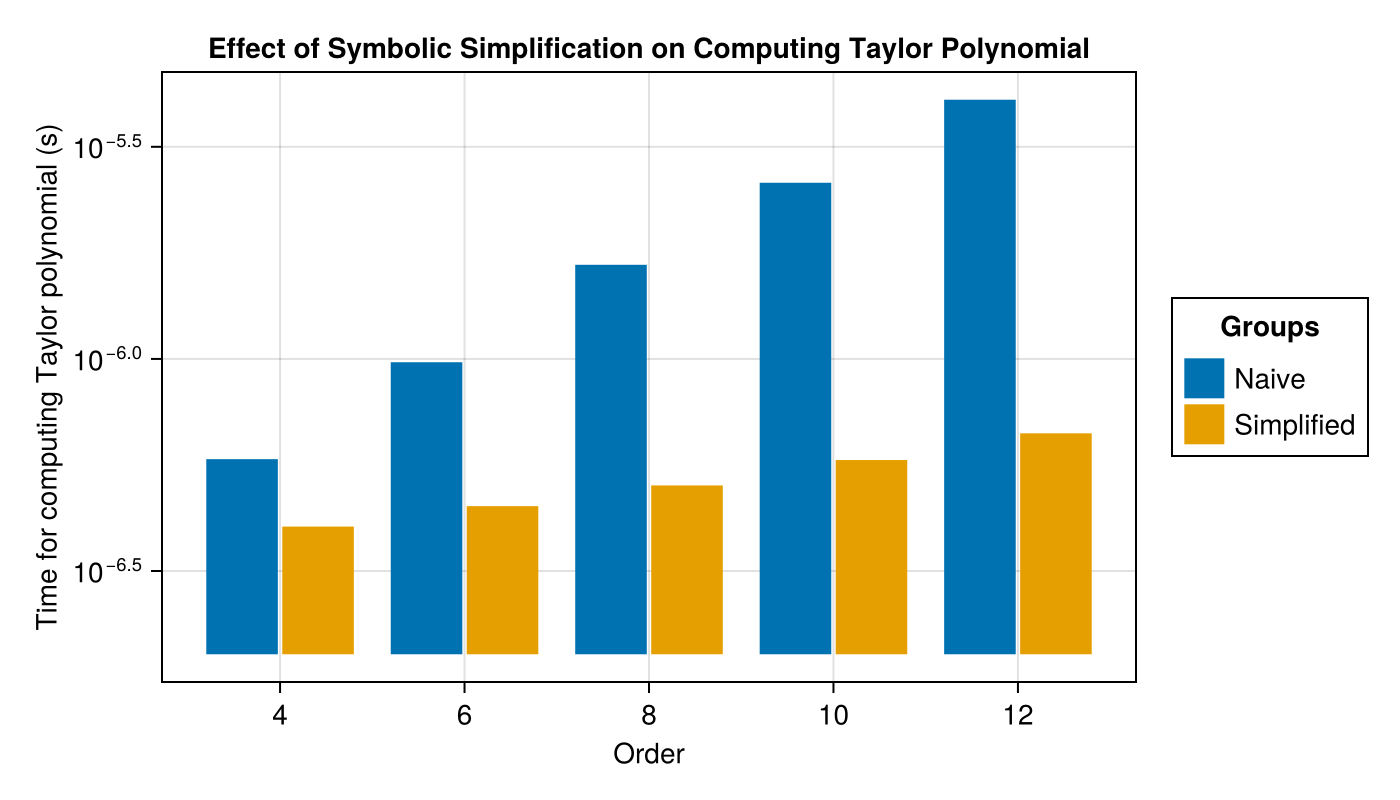}
    \caption{Effect of symbolic compilation on computing Taylor polynomial. For various degrees, performing symbolic simplification in advance (yellow bars) can achieve about one order of magnitude acceleration compared to naively computing (blue bars).}
    \label{fig:simp}
\end{figure}

In our implementation, we execute the process above \textbf{symbolically, at compile time} instead of \textbf{numerically, at run time}. Assuming the target polynomial degree is $p$, we will obtain the final symbolic expression of $u_n^{(1)},\cdots,u_n^{(p)}$ as a directed computational graph with input $u_n$ and $t_n$ only. We will then compile this into a numerical function and use that to advance one step at any time $t$ afterwards. The symbolic compilation will be an one-time cost (usually on the order of microseconds or milliseconds, depending on the complexity of ODE functions), which will be amortized for long integrations where the compiled function will be called thousands or millions of time. \cref{fig:simp} demonstrates the performance of evaluating the Taylor polynomial with and without symbolic compilation, at $t=0$ for the solution of planar circular restricted three-body problem (PCR3BP) ~\cite{murray1999solar}, i.e.
$$
\begin{aligned}
    \dot{x} & = p_x + y,\\
    \dot{y} & = p_y - x,\\
    \dot{p_x} & = - \frac{(1-\mu)(x-\mu)}{((x-\mu)^2+y^2)^{3/2}} - \frac{\mu(x+1-\mu)}{((x+1-\mu)^2+y^2)^{3/2}} + p_y,\\
    \dot{p_y} & = - \frac{(1-\mu)y      }{((x-\mu)^2+y^2)^{3/2}} - \frac{\mu y       }{((x+1-\mu)^2+y^2)^{3/2}} - p_x.
\end{aligned}
$$
where $\mu=0.01$, and $x_0=-0.80,y_0=0.0,p_{x_0}=0.0,p_{y_0}=-0.63$.
Results indicate that symbolic simplification significantly speedups the computation of Taylor polynomials, especially for higher orders.

While we couldn't find an open source implementation of Taylor methods ODE solver based on JAX higher order AD~\footnote{\url{https://docs.jax.dev/en/latest/jax.experimental.jet.html}} to benchmark against, based on the results in their report~\cite{bettencourt_taylor-mode_2022}, we would expect that the coefficient doubling method results in similar performance as the naive implementation without symbolic simplification here.

\subsection{Adaptive polynomial degree and adaptive step size mechanisms}
We first discuss the situation where only step size is adjustable. The Lagrange remainder of $T_n(t)$ could be written as
$$
R_p=\frac{(\mathrm d^{p+1}u/\mathrm dt^{p+1})(\xi)}{(p+1)!}(t-t_n)^{p+1}\approx \frac{u_{n}^{(p+1)}}{(p+1)!}(t-t_n)^{p+1}
$$
where we assume that the first term can be approximated using one additional order of derivative at time $t_n$. Ideally, the remainder $R_p$ should be less than the error tolerance. Therefore, given the current integration timestep $h$, we can define the scaled error $e_p$ based on the absolute error tolerance $T_a$ and the relative error tolerance $T_r$, which is
$$
e_p=\left\|\frac{u_n^{(p+1)}h^{p+1}/(p+1)!}{T_a+|u_n|T_r}\right\|
$$
which should be less than 1 for accurate integration of the ODE system. If $e_p>1$, then this step is rejected and the system is reverted to the previous state. If $e_p\le 1$, it is then processed by a controller $C$, which could be one of Integral Controller (I), Proportional-Integral Controller (PI) or Proportional-Integral-Derivative Controller (PID), which will propose a new step size $h_{\rm new}$ based on current step size $h$ and scaled error $e_p$. Overall, the functionality of controllers could be written as $h_{\rm new}=C(h,e_p)$. This mechanism is summarized in \cref{alg:adapttime}:

\begin{algorithm}
\caption{Adaptive step size integration}
\label{alg:adapttime}
\begin{algorithmic}
\STATE{Define $\dot u=f(u)$, integration time range $t\in(t_0,t_{\rm end})$, initial value $u=u_0$}
\STATE{Let $t=t_0$ and guess an initial step size $h$}
\WHILE{$t<t_{\rm end}$}
\STATE{Integrate one step size $h$ to evolve the solution $u\leftarrow \operatorname{TaylorIntegration}(f,u,h)$}
\STATE{Estimate scaled error $e_p$}
\IF{$e_p>1$}
\STATE{Reject the step size and revert $u$}
\ELSE
\STATE{Update step size $h\leftarrow C(h,e_p)$}
\ENDIF
\ENDWHILE
\RETURN $u$
\end{algorithmic}
\end{algorithm}

We now discuss the situation where both the degree of Taylor polynomial and step size are adjustable, using similar ideas as in~\cite{ekanathan_fully_2024,utkarsh_parallelizing_2022}. This is a special opportunity for Taylor integration method since the degree can be easily adjusted. Assuming the degree $p$ is allowed to vary in the range $[p_{\min},p_{\max}]$, and the operations to compute degree $p$ polynomial is proportional to $p^2$ as a result of using Taylor-mode automatic differentiation. For each value of $p$, the controller will propose a $h(p)$ based on the error tolerance, therefore we can define a work-per-time ratio $r=p^2/h(p)$ to capture the computational cost per unit of time, and find the $p$ value that minimize this ratio. In practice, to avoid large oscillations in degree $p$, after each step the integrator is allowed to re-choose the degree among $\{p-1,p,p+1\}$. This mechanism is summarized in \cref{alg:adapttimeorder}.

\begin{algorithm}
\caption{Adaptive polynomial degree and step size integration}
\label{alg:adapttimeorder}
\begin{algorithmic}
\STATE{Define $\dot u=f(u)$, integration time range $t\in(t_0,t_{\rm end})$, initial value $u=u_0$}
\STATE{Let $t=t_0$ and guess an initial step size $h$}
\WHILE{$t<t_{\rm end}$}
\STATE{Integrate one step size $h$ to evolve the solution $u\leftarrow \operatorname{TaylorIntegration}(f,u,h)$}
\STATE{Initialize the optimal work-per-time ratio $r_{\rm optimal}=+\infty$, optimal degree $p_{\rm optimal}=p$, and the step size corresponding to optimal degree $h_{\rm optimal}$}
\FORALL{$p_{\rm next}\in[\max\{p_{\min} ,p-1\},\min\{p_{\max} ,p+1\}]$}
\STATE{Estimate scaled error $e_p$}
\STATE{Propose step size $h_{\rm next}= C(h,e_p)$}
\STATE{Compute ratio $r=p_{\rm next}^2/h_{\rm next}$}
\IF{$r<r_{\rm optimal}$}
\STATE{Update $p_{\rm optimal}\leftarrow p_{\rm next}$, $r_{\rm optimal}\leftarrow r$}
\ENDIF
\ENDFOR
\IF{$e_p>1$}
\STATE{Reject the step size and revert $u$}
\ELSE
\STATE{Update step size $h\leftarrow h_{\rm optimal}$, degree $p\leftarrow p_{\rm optimal}$}
\ENDIF
\ENDWHILE
\RETURN $u$
\end{algorithmic}
\end{algorithm}

\section{Experiments}
\label{sec:experiments}
\subsection{Solving various non-stiff ODE problems}
We focus on non-stiff problems, where explicit Taylor methods are theoretically competitive. We test several Taylor methods with different fixed polynomial degree (6, 8, 10, 12) against several Runge-Kutta methods (DP5\cite{dormand1980family}, Tsit5\cite{tsitouras2011runge}, Vern6, Vern8~\cite{verner2010numerically}), on the following four non-stiff ODE problems:
\paragraph{Lotka-Volterra problem}
$$
\begin{cases}
\dot x= ax - bxy\\
\dot y= -cy + dxy
\end{cases}
$$
where $t\in(0.0, 10.0)$, $a=1.5, b=1.0,c=3.0,d=1.0$, and $x_0=y_0=1.0$.
\paragraph{Fitzhugh-Nagumo problem}
$$
\begin{cases}
\dot v = v - v^3/3 -w + l\\
\dot w = (v +  a - bw)/\tau
\end{cases}
$$
where $t\in(0.0, 10.0)$, $a=0.7,b=0.8,\tau=12.5,l=0.5$ and $v_0=w_0=1.0$.
\paragraph{Rigid body problem}
$$
\begin{cases}
\dot x  = I_1yz\\
\dot y  = I_2zx\\
\dot z  = I_3xy + \sin(t)^2/4
\end{cases}
$$
where $t\in(0.0, 10.0)$, $I_1=-2.0,I_2=1.25,I_3=-0.5$ and $x_0=1.0,y_0=0.0,z_0=0.9$.
\paragraph{Random linear problem}
$$
\dot u=Au
$$
where $t\in(0.0, 10.0)$, $A\in\mathbb R^{16\times 16}$ and $u_0\in\mathbb R^{16}$, both of their elements are sampled from normal distribution $N(0,1)$.
\begin{figure}[htbp]
    \centering
    \begin{subfigure}{0.49\textwidth}
        \centering
        \includegraphics[width=\textwidth]{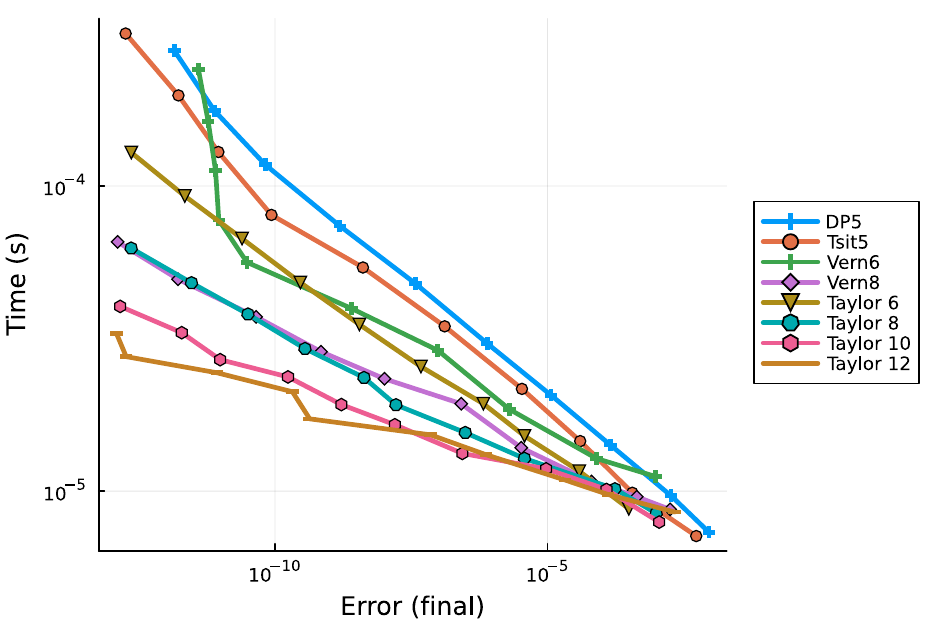}
        \caption{Lotka-Volterra problem}
        \label{fig:sub1}
    \end{subfigure}
    \hfill
    \begin{subfigure}{0.49\textwidth}
        \centering
        \includegraphics[width=\textwidth]{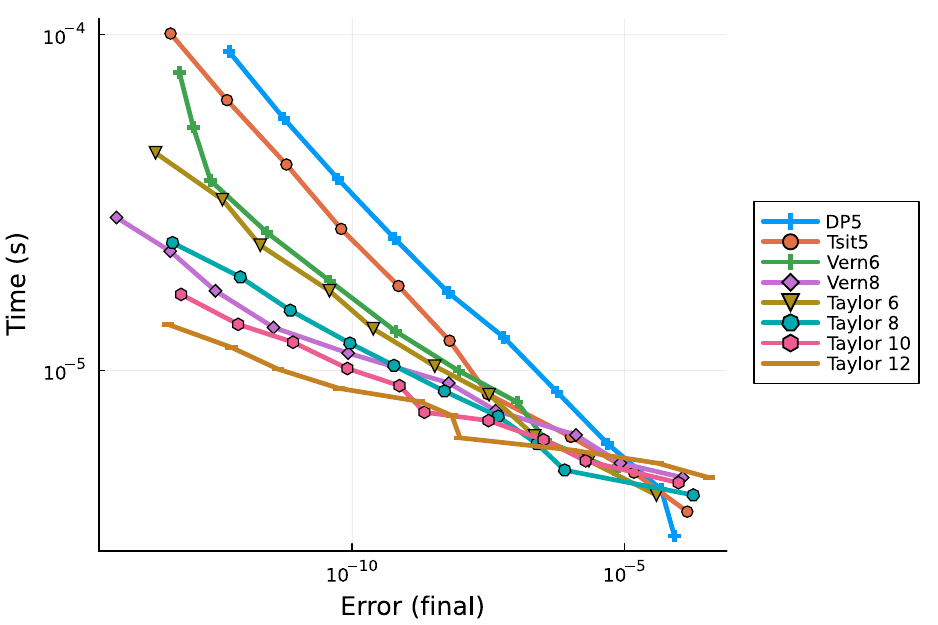}
        \caption{Fitzhugh-Nagumo problem}
        \label{fig:sub2}
    \end{subfigure}

    \vspace{0.5em} % 调整上下排间距（可选）

    \begin{subfigure}{0.49\textwidth}
        \centering
        \includegraphics[width=\textwidth]{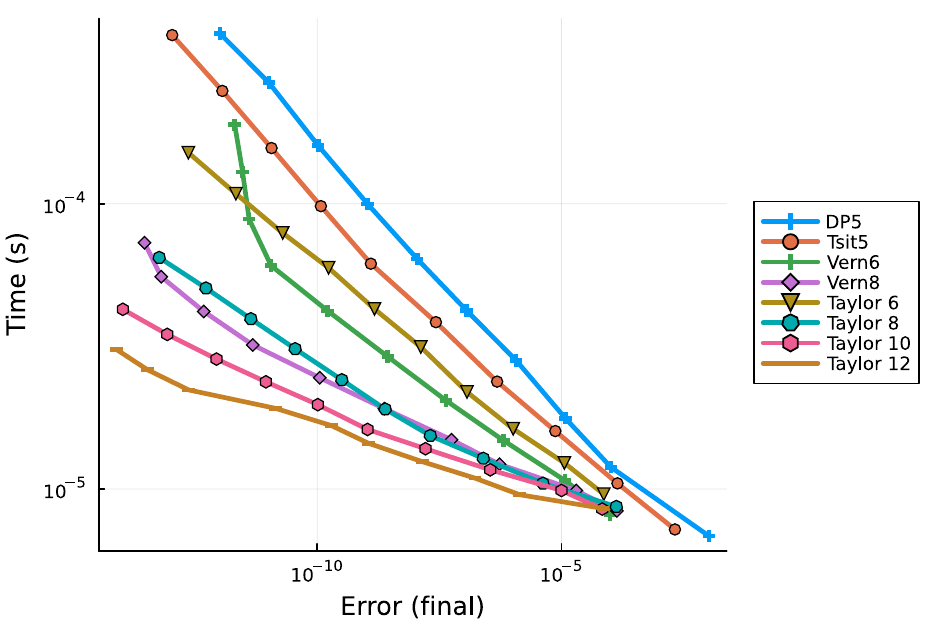}
        \caption{Rigid body problem}
        \label{fig:sub3}
    \end{subfigure}
    \hfill
    \begin{subfigure}{0.49\textwidth}
        \centering
        \includegraphics[width=\textwidth]{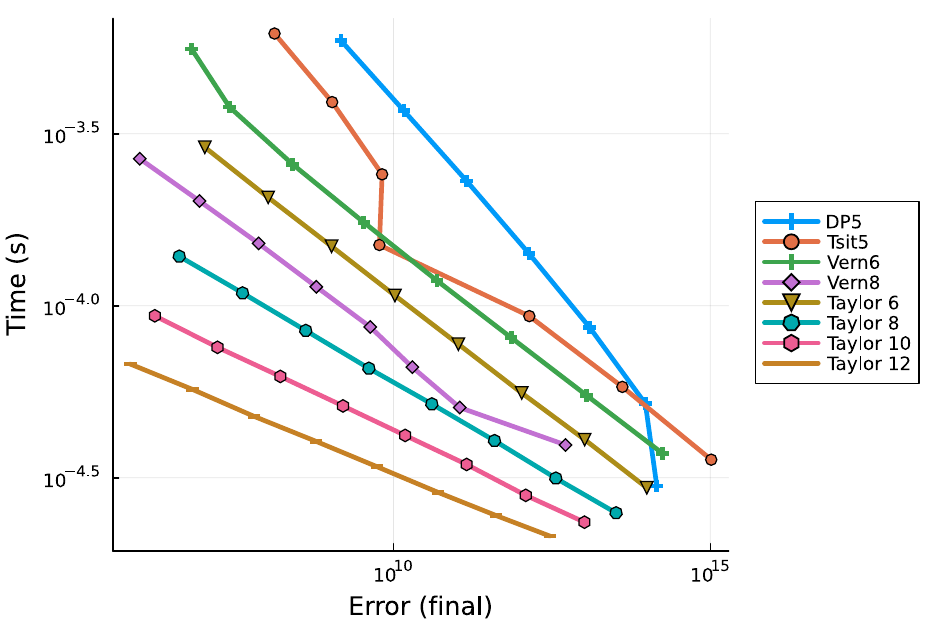}
        \caption{Random linear problem}
        \label{fig:sub4}
    \end{subfigure}
    
    \caption{Work-precision diagrams for Taylor methods and Runge-Kutta methods in various non-stiff ODE problems}
    \label{fig:nonstiff}
\end{figure}

While the work-precision curve is similar for same order Taylor methods and Runge-Kutta methods (such as Taylor6 vs Vern6, Taylor8 vs Vern8), the degree of Taylor polynomial (and therefore the convergence order) of Taylor methods can be systematically increased to achieve higher accuracy. Therefore, higher order Taylor methods can have advantage over Runge-Kutta methods, especially at higher precision.

\subsection{Adaptive polynomial degree tests}

We consider the following problem:
$$
\dot y=\phi'(t)(y-y^2)
$$
where
$$
\phi(t)=\kappa\tanh\frac{t-t_0}\delta
$$
and $c=2.0,\delta=1.0,t_0=5.0,\kappa=10.0$ and  $y_0=1/(1+c)$. The analytical solution will be
$$
y(t)=\frac1{1+c\exp\left(-\phi(t)\right)}
$$
which has singularity at $\phi(t)=\ln c\pm \pi i$. Therefore, for some $t$ on the real axis, the convergence radius of Taylor series will be
$$
\rho(t)\approx\frac{\pi}{|\phi'(t)|}
$$
Now, since
$$
\phi'(t)=
\frac\kappa\delta \cosh^{-2}\left(\frac{t-t_0}{\delta}\right)
$$
has a peak at $t_0$, it is expected to have small convergence radius near that (the ``hard'' region) and large convergence radius away from that (the ``easy'' region). In the hard region, using a large polynomial degree $p$ cannot increase the step size $h$ since it is controlled by the convergence radius, and therefore smaller $p$ is preferred; in the easy region, using a large $p$ can greatly increase the step size, and therefore larger $p$ is preferred. In the work-precision diagram below, we use four different fixed-degree Taylor methods (degree 6, 8, 10 and 12 respectively) and an adaptive-degree Taylor method (degree can vary between 6 and 12) to demonstrate this idea.

\begin{figure}[!ht]
    \centering
    \includegraphics[width=0.7\linewidth]{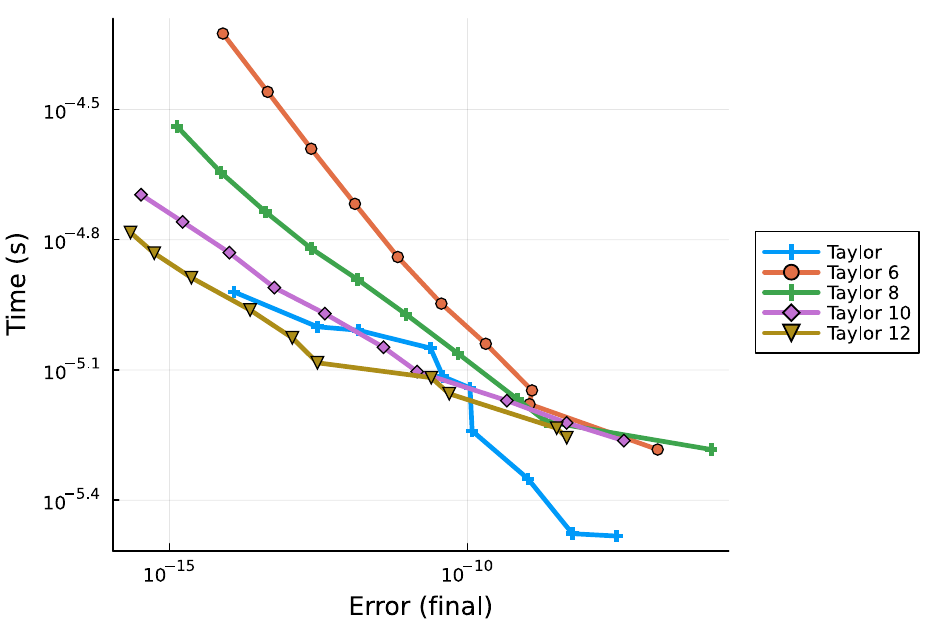}
    \caption{Comparison of adaptive-degree Taylor methods (degree can vary between 6 and 12) and fixed-degree Taylor method (6, 8, 10 and 12). At lower tolerances, the adaptive-degree method has a chance to use different degree at different region to achieve better performance.}
    \label{fig:adaptive}
\end{figure}

\section{Conclusions}
\label{sec:conclusions}

In this paper, we introduced explicit Taylor methods and conducted a comprehensive review of existing implementations of it. Pointing out the problems with existing implementations, we proposed a new implementation which leveraged a general purpose higher-order automatic differentiation engine and a computer algebra system, making it efficient, broadly applicable and user-friendly. We then introduced the mechanisms to adjust Taylor polynomial degree and step size at different regions across the integration time period, which help the method to achieve better performance. In numerical experiments, we compared the fixed-degree Taylor methods with various Runge-Kutta methods on four different non-stiff ODE problems. The work-precision diagrams demonstrated the advantage of Taylor methods, especially at higher precision since the convergence order of Taylor methods could be arbitrarily high. We then proposed a problem where very different convergence radii of Taylor series are expected at different regions, making the mechanism of adaptive polynomial degree profitable. In this case, the adaptive-degree Taylor method outperforms the fixed-degree ones.

However, as the same for all explicit methods, explicit Taylor methods are generally unsuitable for problems that are significantly stiff, since it might be forced to use very small time steps even with a high degree of Taylor polynomial. We hope to extend our approach toward a better implementation of implicit Taylor methods~\cite{scott_solving_2000} as well.

\section*{Acknowledgments}
This material is based upon work supported by the U.S. National Science Foundation under award Nos CNS-2346520,  RISE-2425761, and DMS-2325184, by the Defense Advanced Research Projects Agency (DARPA) under Agreement No. HR00112490488,  by the Department of Energy, National Nuclear Security Administration under Award Number DE-NA0003965, by the U.S. Department of Energy, Office of Science, Office of Advanced Scientific Computing Research as part of the Next Generation of Scientific Software Technologies program, Stewardship and Advancement of Programming Systems and Tools (S4PST) project, under contract number CW54751 and by the United States Air Force Research Laboratory under Cooperative Agreement Number FA8750-19-2-1000.  Neither the United States Government nor any agency thereof, nor any of their employees, makes any warranty, express or implied, or assumes any legal liability or responsibility for the accuracy, completeness, or usefulness of any information, apparatus, product, or process disclosed, or represents that its use would not infringe privately owned rights. Reference herein to any specific commercial product, process, or service by trade name, trademark, manufacturer, or otherwise does not necessarily constitute or imply its endorsement, recommendation, or favoring by the United States Government or any agency thereof. The views and opinions of authors expressed herein do not necessarily state or reflect those of the United States Government or any agency thereof." The views and conclusions contained in this document are those of the authors and should not be interpreted as representing the official policies, either expressed or implied, of the United States Air Force or the U.S. Government.

\appendix
\section{Comparison with TaylorIntegration.jl}
\begin{figure}[!ht]
    \centering
    \includegraphics[width=0.7\linewidth]{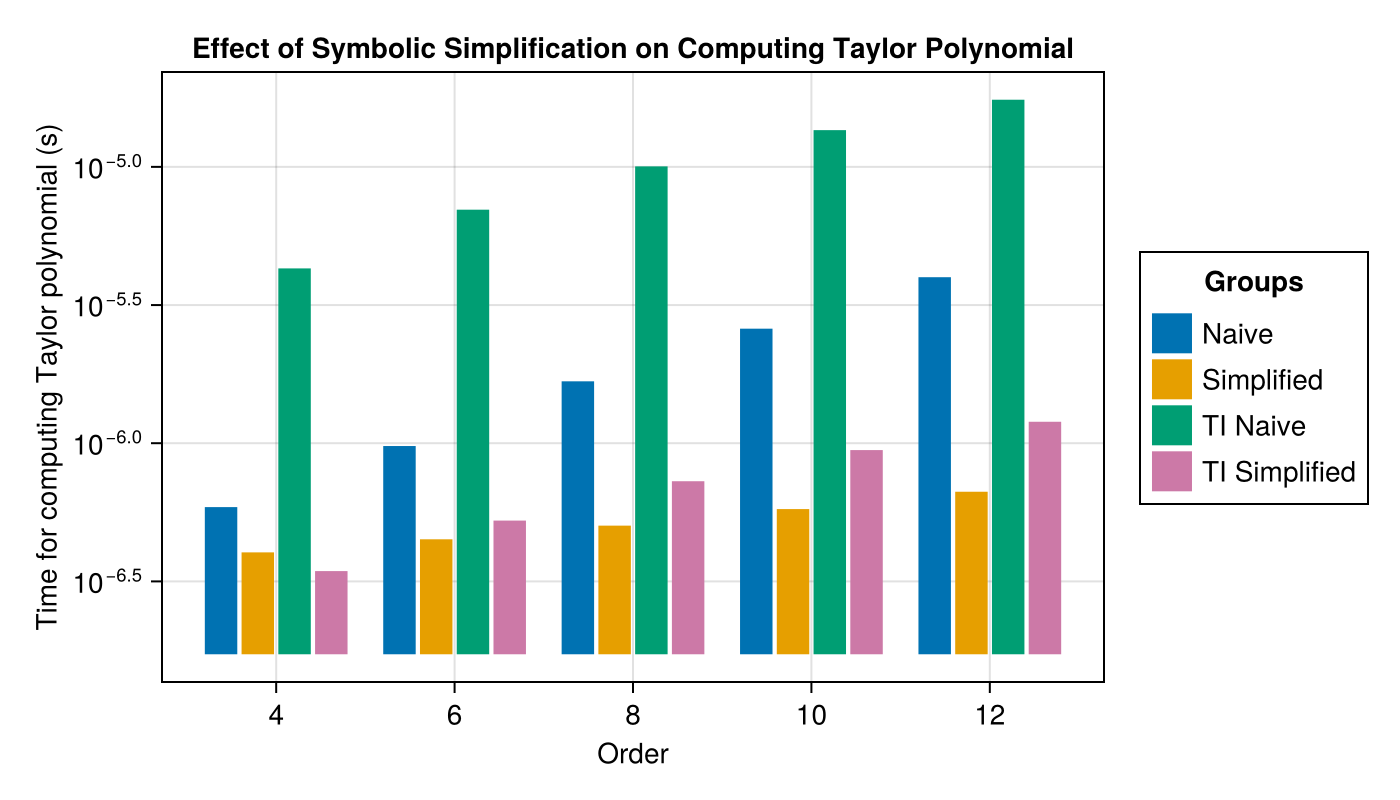}
    \caption{Effect of symbolic compilation on computing Taylor polynomial, with both this work and TaylorIntegration.jl. In addition to this work (where naive and simplified are shown in blue and yellow bars), computing with TaylorIntegration.jl is shown in green bars for naive and pink for simplified.}
    \label{fig:simpti}
\end{figure}

Here we present more detail about the comparison with TaylorIntegration.jl, which is another implementation of Taylor methods in Julia that we mentioned before. Based on the same PCR3BP problem, in \cref{fig:simpti} our implementation shows better performance both before and after simplification, which demonstrates that a general and systematic computer algebra system can perform better than \emph{ad hoc} code transformers. In addition, while TaylorIntegration.jl requires user to use macros to transform the ODE function as well as rewriting the function in a compatible way, our approach removes the need for any user-facing code modification, and should be readily applied to existing codebases.

In addition to integration with floating point numbers, we note that TaylorIntegration.jl also supports interval arithmetic so it could have other purposes like global error analysis.
\nocite{*}
\bibliographystyle{siamplain}
\bibliography{references}
\end{document}